\documentclass[12pt,letter]{amsart}
\usepackage[latin1]{inputenc}
\usepackage{amsmath,amsthm,amssymb,latexsym}
\usepackage[T1]{fontenc}
\usepackage{libertine} 
\usepackage{euler} 
\usepackage{eucal} 
\usepackage{color}
\usepackage{graphicx}
\usepackage{pb-diagram}
\usepackage{pst-3d,pst-coil,pst-eps,pst-fill,pst-grad,pst-math,pst-text,pst-tree}

\newcommand\C[1]{\mathcal{#1}}

\newtheorem{proposition}{Proposition}[section]
\newtheorem{theorem}[proposition]{Theorem}

\newtheorem{assumption}[proposition]{Assumption}

\newtheorem{corollary}[proposition]{Corollary}
\newtheorem{fact}[proposition]{Fact}
\newtheorem{claim}[proposition]{Claim}

\theoremstyle{definition}
\newtheorem{definition}[proposition]{Definition}

\newtheorem{remark}[proposition]{Remark}

\def\monster{\mathbb{M}}
\def\E{\varepsilon}
\def\dist{\mathbf{d}}
\newcommand{\indep}{\makebox[13pt]{\raisebox{-0.75ex}{$\smile$}\hspace{-3.50mm}
$\mid$\, }}

\def\ordencea{\prec_{\C{K}}}
\def\ordenceac{\succ_{\C{K}}}
\newcommand{\eop}[1]{
\hspace{10mm} \vspace{-6mm}
\begin{flushright}
\qedsymbol$_{\text{#1}}$\\ \ \\
\end{flushright}
}
\newenvironment{prueba}[1][{\it Proof}]{\noindent {\it #1.} }{}
\newcommand{\bdem}[1][Proof]{\begin{prueba}[#1]}
\newcommand{\edem}[1][]{\eop{#1}
\end{prueba}}
\def\bsdem{\begin{prueba}[Reference]}\def\bsindem{\begin{proof}[\ ]}

\def\tuple{\overline}
\def\rest{\upharpoonright}

\mathchardef\mhyphen="2D
\newcommand{\gatp}{{\sf ga\mhyphen tp}}
\newcommand{\gaS}{{\sf ga\mhyphen S}}

\begin{document}
\title[A stability transfer theorem in d-Tame Metric AEC $\cdots$]{A stability transfer theorem in d-tame metric abstract elementary classes}
\author[P. Zambrano]{Pedro Zambrano}
\address{{\rm E-mail:} {\it phzambranor@unal.edu.co}\\
Departamento de Matem\'aticas Universidad Nacional de Colombia, AK
30 $\#$ 45-03, Bogot\'a - Colombia}
\thanks{
The author wants to thank Andr\'es Villaveces for his suggestions in early proofs of the main theorem of this paper. The author is very thankfully to Tapani Hyttinen for the nice discussions and suggestions about this paper during a short visit to Helsinki in 2009. The author
  was partially supported by Colciencias.}
\date{August 1st, 2011}
\begin{abstract}
In this paper, we study a stability transfer theorem in $d$-tame Metric Abstract Elementary classes, in a similar way as in \cite{BaKuVa}, but using superstability-like assumptions which involves a new independence notion ({\it Tame Independence}) instead of $\aleph_0$-locality.
\end{abstract}
\maketitle

\section{Introduction}

Discrete {\it tame} Abstract Elementary Classes are a very special kind of Abstract Elementary Classes (shortly, AECs) which have a categoricity transfer theorem (see \cite{GrVa2}) and a nice stability transfer theorem (see \cite {BaKuVa}). In fact -under $\aleph_0$-tameness and $\aleph_0$-locality (assuming $LS(\C{K})=\aleph_0$)-, J. Baldwin, D. Kueker and M. VanDieren proved in  \cite{BaKuVa} that $\aleph_0$-Galois-stability implies $\kappa$-Galois-stability for every cardinality $\kappa$. First, they proved that $\aleph_0$-Galois-stability implies $\aleph_n$-Galois stability for every $n<\omega$ (in fact, their argument works for getting $\kappa$-Galois-stability if $cf(\kappa)>\omega$) and so (by $\aleph_0$-locality) $\aleph_\omega$-Galois-stable (where the same argument works for getting $\kappa$-Galois stability if $cf(\kappa)=\omega$).
\\ \\
\indent {\it Metric Abstract Elementary Classes} (for short, MAECs) correspond to a kind of amalgam between {\it AECs} and {\it Continuous Logic Elementary Classes}, although we drop uniformly continuity of the symbols of the languages (for our purposes, it is enough to take closed functions). In this setting, it is enough to consider dense subsets of the models, so this is the reason because all our analysis considers density character instead of cardinality of the models. In general, we can define a distance between Galois-types in this setting, which is a metric under suitable assumptions (see \cite{Hi,ViZa}). Because of that, we adapt a notion of {\it Tameness} using these new tools given in this setting.
\\ \\
\indent In section 2, we study a suitable notion of independence (which we call {\it Tame Independence}) which we use for proving the stability transfer theorem in this setting. This is one of the differences between our paper and \cite{BaKuVa} -they just used a combinatoric argument to get their result-.  in this paper, also we strongly use superstability-like assumptions ($\E$-locality, assumption \ref{superstability_tameness}) to get our main theorem.
\\ \\
In section 3, we provide the proof of our main result of stability transfer theorem, which -roughly speaking- says that under $d$-tameness, $\aleph_0$ and $\aleph_1$-d-stability and some suitable superstablity-like assumptions -via tame independence- we have $\kappa$-d-stability for all cardinality $\kappa$.

\section[Independence in $d$-tame MAEC]{An independence notion in $d$-tame metric abstract elementary classes.}\label{tameMAEC}
In this section, we provide a definition of tameness adapted to the setting of metric abstract elementary classes and a suitable notion of independence, which we will use in section 3 for proving an upward stability transfer theorem.

This section is devoted to develop a suitable notion of stability towards proving the following fact:
\\ \\
{\bf Theorem \ref{stab_transfer2}.}
{\it Let $\C{K}$ be a 
$\mu$-$\mathbf{d}$-tame (for some $\mu<\kappa$) MAEC. Suppose that 
$\C{K}$ is $[LS(\C{K}),\kappa)$-cofinally-d-stable.
Define $$\lambda:=\min\{\theta<\kappa: \mu< \theta \text{ and  $\C{K}$ is $\theta$-$\mathbf{d}$-stable }\},$$  $$\zeta:=\min\{\xi: 2^{\xi}>\lambda\}$$ and $$\zeta^*:=\max\{\mu^+,\zeta\}.$$ If $cf(\kappa)\ge \zeta^*$ then $\C{K}$ is $\kappa$-$\mathbf{d}$-stable.}
\\ \\
\indent We will provide a proof of theorem \ref{stab_transfer2} in section 3.
\\ \\
\indent Under superstability-like assumptions ($\E$-locality) on a notion of independence which we will define in this section, the theorem above implies $\kappa$-d-stability for every $\kappa$. 

For the basic notions and facts in MAECs, we refer the reader to \cite{Hi,ViZa}. For the sake of completeness, we provide some of the most relevant notions and facts which we use in this paper.

\begin{definition}
Let $(X,\tau)$ be a topological space. The {\it density character} of $(X,\tau)$ is defined as the minimum cardinality of a dense subset of $X$.
\end{definition} 

\begin{definition}[distance between Galois types]
Let $\C{K}$ be an MAEC with AP and JEP -so Galois types over a model $M$ correspond to orbits of automorphisms of a fixed monster model $\monster$ which fix $M$ pointwise-. Let $M\in \C{K}$ and $p,q\in \gaS(M)$. Define $d(p,q):=\inf\{d(a,b): a,b\in\monster, a\models p \text{\ and } b\models q \}$.
\end{definition}

\begin{definition}
Let $\C{K}$ be an MAEC with AP and JEP. We say that $\C{K}$ has the {\it Continuity Type Property}\footnote{CTP is called {\it Perturbation Property} in \cite{Hi}} (for short, CTP) iff for any convergent sequence $(a_n)_{n<\omega}$ in $\monster$, if $(a_n)\to a$ and $\gatp(a_n/M)=\gatp(a_0/M)$ for all $n<\omega$, then $\gatp(a/M)=\gatp(a_0/M)$.
\end{definition}

\begin{fact}[Hirvonen-Hyttinen]
Let $\C{K}$ be an MAEC with AP and JEP. $d$ defined as above is a metric iff $\C{K}$ has the CTP.
\end{fact}

Most of the natural examples (e.g., Banach Spaces and Elementary Continuous Logic Classes) satisfy CTP. So, we may assume that distance between Galois types is in fact a metric.

\begin{definition}[$\lambda$-d-stability]
Let $\C{K}$ be an MAEC with AP and JEP and $\lambda\ge LS(\C{K})$. We say that $\C{K}$ is {\it $\lambda$-d-stable} iff given any $M\in \C{K}$ with density character $\lambda$, $dc(\gaS(M))\le \lambda$
\end{definition}

\begin{definition}[Cofinal-d-stability]
Let $\C{K}$ be an MAEC with AP and JEP and $LS(\C{K})\le \lambda <\kappa$. We say that $\C{K}$ is {\it $[\lambda,\kappa)$-cofinally-d-stable} iff given $\theta\in [\lambda,\kappa)$ there exists $\theta'\ge \theta$ in $[\lambda,\kappa)$ such that $\C{K}$ is $\theta'$-d-stable.
\end{definition}

\begin{definition}[Universality]
Let $\C{K}$ be an MAEC and $M\ordencea N$ in $\C{K}$. We say that $N$ is {\it $\mu$-d-universal over $M$} iff  for every $M'\ordenceac M$ of density character $\mu$ there exists a $\C{K}$-embedding $f:M'\to N$ which fixes $M$ pointwise. We say that $N$ is d-universal over $M$ iff it is $dc(M)$-d-universal. We drop d if the metric context is clear.
\end{definition}

Under d-stability, universal models exist. 

\begin{fact}\label{Existence_Universal}
Let $\C{K}$ be a $\mu$-d-stable MAEC. Given $M\in \C{K}$ of density character $\mu$, there exists $M'\ordenceac M$ universal over $M$.
\end{fact}

$\mu$-Tameness in (discrete) AECs says that the difference between two Galois-types $p,q\in \gaS(M)$ is given by some $N\ordencea M$ of size $\mu$. Since in this setting we have a distance between Galois-types (see \cite{Hi}), so we adapt this notion to the metric setting. 

\begin{definition}[$\mathbf{d}$-tameness]\label{tameness}
Let $\C{K}$ be a MAEC and $\mu\ge LS(\C{K})$. We say that $\C{K}$ is {\it $\mu$-$\mathbf{d}$-tame}
iff for every $\E>0$ there exists $\delta_\E>0$ such that if for any $M\in \C{K}$ of density character $\ge\mu$ we have that $\dist(p,q)\ge \E$ where $p,q\in \gaS(M)$, then there exists $N\ordencea M$ of density character $\mu$ such that $\dist(p\rest N,q\rest N)\ge \delta_\E$.
\end{definition}

\begin{assumption}
The definitions given below use $\lambda$, $\mu$ and $\zeta^*$ defined above. So, throughout this section, we assume that $\C{K}$ is a $\mu$-d-tame and a $\lambda$-d-stable MAEC. Also, we suppose that $\C{K}$ satisfies AP and JEP, so we may able to construct a homogeneous monster model $\monster\in \C{K}$ and we consider the Galois-types over $M\in \C{K}$ as orbits under  $Aut(\monster/M)$.
\end{assumption}

As we did in the definition of $d$-tameness, we can adapt the notion of splitting to MAECs using the distance between Galois-types.

\begin{definition}
Let $N\ordencea M$ and $\E>0$. We say that $\gatp(a/M)$
{\it tame-$<\zeta^*$-$\E$-splits} over $N$ iff for every submodel $N'\ordencea N$ with density character $<\zeta^*$, there are models $N'\ordencea N_1,N_2\ordencea
M$ with density character $<\zeta^*$ and $h:N_1\cong_{N'} N_2$ such that
$\dist(\gatp(a/N_2),h(\gatp(a/N_1))\ge \E$. If it is clear, we drop $<\zeta^*$ and we just say that $\gatp(a/M)$ tame-$\E$-splits over $N$.
If $\gatp(a/M)$ does not tame-$\varepsilon$-split over $N$, we denote that by $a\indep^{T,\E}_N
M$.
\end{definition}

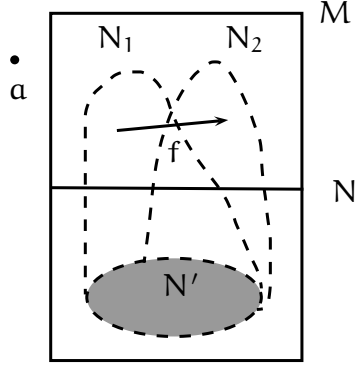
\begin{center}
\scalebox{1} 
{
\begin{pspicture}(0,-2.3992188)(4.9228125,2.4392188)
\definecolor{color7b}{rgb}{0.6,0.6,0.6}
\psframe[linewidth=0.04,dimen=outer](4.0209374,2.2607813)(0.6409375,-2.3992188)
\psline[linewidth=0.04cm](0.6609375,-0.07921875)(4.0209374,-0.09921875)
\psellipse[linewidth=0.04,linestyle=dashed,dash=0.16cm 0.16cm,dimen=outer,fillstyle=solid,fillcolor=color7b](2.3009374,-1.5392188)(1.18,0.54)
\pscustom[linewidth=0.04,linestyle=dashed,dash=0.16cm 0.16cm]
{
\newpath
\moveto(1.1209375,-1.4992187)
\lineto(1.1209375,-0.41921875)
\curveto(1.1209375,0.12078125)(1.1209375,0.7657812)(1.1209375,0.87078124)
\curveto(1.1209375,0.97578126)(1.1959375,1.1457813)(1.2709374,1.2107812)
\curveto(1.3459375,1.2757813)(1.4559375,1.3707813)(1.4909375,1.4007813)
\curveto(1.5259376,1.4307812)(1.6059375,1.4607812)(1.6509376,1.4607812)
\curveto(1.6959375,1.4607812)(1.7959375,1.4607812)(1.8509375,1.4607812)
\curveto(1.9059376,1.4607812)(2.0209374,1.3857813)(2.0809374,1.3107812)
\curveto(2.1409376,1.2357812)(2.2159376,1.1007812)(2.2309375,1.0407813)
\curveto(2.2459376,0.98078126)(2.3159375,0.81078124)(2.3709376,0.7007812)
\curveto(2.4259374,0.5907813)(2.5459375,0.39078125)(2.6109376,0.30078125)
\curveto(2.6759374,0.21078125)(2.7959375,0.05078125)(2.8509376,-0.01921875)
\curveto(2.9059374,-0.08921875)(2.9859376,-0.21421875)(3.0109375,-0.26921874)
\curveto(3.0359375,-0.32421875)(3.0709374,-0.40421876)(3.0809374,-0.42921874)
\curveto(3.0909376,-0.45421875)(3.1059375,-0.49921876)(3.1109376,-0.51921874)
\curveto(3.1159375,-0.5392187)(3.1409376,-0.5992187)(3.1609375,-0.63921875)
\curveto(3.1809375,-0.67921877)(3.2259376,-0.77921873)(3.2509375,-0.83921874)
\curveto(3.2759376,-0.89921874)(3.3109374,-0.9792187)(3.3209374,-0.99921876)
\curveto(3.3309374,-1.0192188)(3.3609376,-1.0792187)(3.3809376,-1.1192187)
\curveto(3.4009376,-1.1592188)(3.4259374,-1.2242187)(3.4309375,-1.2492187)
\curveto(3.4359374,-1.2742188)(3.4459374,-1.3192188)(3.4509375,-1.3392187)
\curveto(3.4559374,-1.3592187)(3.4609375,-1.3892188)(3.4609375,-1.4192188)
}
\pscustom[linewidth=0.04,linestyle=dashed,dash=0.16cm 0.16cm]
{
\newpath
\moveto(1.9009376,-1.0792187)
\lineto(1.9409375,-0.7892187)
\curveto(1.9609375,-0.64421874)(1.9909375,-0.37921876)(2.0009375,-0.25921875)
\curveto(2.0109375,-0.13921875)(2.0409374,0.09578125)(2.0609374,0.21078125)
\curveto(2.0809374,0.32578126)(2.1209376,0.50578123)(2.1409376,0.57078123)
\curveto(2.1609375,0.6357812)(2.2109375,0.78078127)(2.2409375,0.86078125)
\curveto(2.2709374,0.94078124)(2.3559375,1.1107812)(2.4109375,1.2007812)
\curveto(2.4659376,1.2907813)(2.5559375,1.4107813)(2.5909376,1.4407812)
\curveto(2.6259375,1.4707812)(2.6859374,1.5207813)(2.7109375,1.5407813)
\curveto(2.7359376,1.5607812)(2.7809374,1.5857812)(2.8009374,1.5907812)
\curveto(2.8209374,1.5957812)(2.8659375,1.5957812)(2.8909376,1.5907812)
\curveto(2.9159374,1.5857812)(2.9809375,1.5307813)(3.0209374,1.4807812)
\curveto(3.0609374,1.4307812)(3.1509376,1.2757813)(3.2009375,1.1707813)
\curveto(3.2509375,1.0657812)(3.3209374,0.9007813)(3.3409376,0.8407813)
\curveto(3.3609376,0.78078127)(3.3909376,0.69578123)(3.4009376,0.67078125)
\curveto(3.4109375,0.6457813)(3.4309375,0.5857813)(3.4409375,0.55078125)
\curveto(3.4509375,0.5157812)(3.4609375,0.43578124)(3.4609375,0.39078125)
\curveto(3.4609375,0.34578124)(3.4659376,0.23578125)(3.4709375,0.17078125)
\curveto(3.4759376,0.10578125)(3.4859376,0.0)(3.4909375,-0.03921875)
\curveto(3.4959376,-0.07921875)(3.5059376,-0.14921875)(3.5109375,-0.17921875)
\curveto(3.5159376,-0.20921876)(3.5259376,-0.29921874)(3.5309374,-0.35921875)
\curveto(3.5359375,-0.41921875)(3.5509374,-0.5342187)(3.5609374,-0.58921874)
\curveto(3.5709374,-0.64421874)(3.5809374,-0.74421877)(3.5809374,-0.7892187)
\curveto(3.5809374,-0.83421874)(3.5809374,-0.92921877)(3.5809374,-0.9792187)
\curveto(3.5809374,-1.0292188)(3.5809374,-1.1292187)(3.5809374,-1.1792188)
\curveto(3.5809374,-1.2292187)(3.5809374,-1.3192188)(3.5809374,-1.3592187)
\curveto(3.5809374,-1.3992188)(3.5709374,-1.4992187)(3.5609374,-1.5592188)
\curveto(3.5509374,-1.6192187)(3.5159376,-1.6842188)(3.4909375,-1.6892188)
\curveto(3.4659376,-1.6942188)(3.4309375,-1.6992188)(3.4009376,-1.6992188)
}
\psline[linewidth=0.04cm,fillcolor=color7b,arrowsize=0.05291667cm 2.0,arrowlength=1.4,arrowinset=0.4]{->}(1.5409375,0.68078125)(3.0209374,0.82078123)
\rput(1.5323437,1.9){$N_1$}
\rput(3.2323437,1.9){$N_2$}
\rput(2.3,0.41078126){$f$}
\rput(2.390625,-1.3692187){$N'$}
\rput(4.572344,-0.12921876){$N$}
\rput(4.4423437,2.2507813){$M$}
\psdots[dotsize=0.12](0.1809375,1.6207813)
\rput(0.22234374,1.2){$a$}
\end{pspicture}
}
\end{center}

\begin{definition}
Let $N\ordencea M$. We say that $a$ is {\it tame-independent} from $M$ over $\C{N}$ iff for every $\E>0$ we have that $a\indep^{T,\E}_N M$. We denote this by $a\indep^T_{N} M$
\end{definition}

In the rest of this section we will prove some basic properties of {\it tame independence}.

\begin{proposition}[Monotonicity]\label{monotonicityt}
Let $M_0\ordencea M_1\ordencea M_2\ordencea M_3$ and suppose that $a\indep^T_{M_0}{M_3}$. Then $a\indep^T_{M_1} M_2$.
\end{proposition}
\bdem
Since $a\indep^T_{M_0}{M_3}$, given $\E>0$ there exists a model $N'\ordencea M_0$ with density character $<\zeta^*$ such that for every models $N'\ordencea N_1\stackrel{h}{\cong}_{N'} N_2\ordencea M_3$ with density character $<\zeta^*$ we have that \linebreak
$d(\gatp(a/N_2),\gatp(h(a)/N_2))<\E$. But we have that $N'\ordencea M_1$ and also it holds in particular if $N'\ordencea N_1\stackrel{h}{\cong}_{N'}N_2\ordencea M_2$. Therefore, $a\indep^T_{M_1} M_2$.
\edem[Prop. \ref{monotonicityt}]

\begin{fact}[Invariance]
Let $f\in Aut(\monster)$. If $a\indep^{T,\E}_N M$ then $f(a)\indep^{T,\E}_{f(N)} f(M)$.
\end{fact}


The following fact strongly uses the $\lambda$-d-stability hypothesis.

\begin{proposition}[Locality]\label{Locality2}
For every $N$, $a$ and every $\E>0$ there exists $M\ordencea N$ of density character $<\zeta^*$ such that $a\indep^{T,\E}_{M} N$.
\end{proposition}
\bdem Suppose that there exists $p:=\gatp(\overline{a}/N)$ such
that $p\not\hspace{-2.5mm}\indep^{T,\E}_M N$ for every $M\ordencea N$
with density character $<\zeta^*$. If $\overline{a}\in N$, it is
straightforward to see that $p$ does not $\varepsilon$-split over
its domain. Then, suppose that $\overline{a}\notin N$.
\\ \\
\indent We will construct a sequence of models $\langle M_\alpha,
N_{\alpha,1}, N_{\alpha,2} : \alpha<\zeta \rangle$ 
in the following
way: First, take $M_0\ordencea N$ as any submodel of density
character $<\zeta^*$.
\\ \\
\indent Suposse $\alpha:=\gamma+1$ and that $M_\gamma$ (with
density character $<\zeta^*$) has been constructed. Therefore $p$
$\varepsilon$-splits over $M_\gamma$. Then there exist
$M_\gamma\ordencea N_{\gamma,1},N_{\gamma,2}\ordencea N $ with
density character $<\zeta^*$ and $F_\gamma:
N_{\gamma,1}\cong_{M_\gamma} N_{\gamma,2}$ such that
$d(F_\gamma(p\rest N_{\gamma,1}),p\rest N_{\gamma,2})\ge
\varepsilon$. Take $M_{\gamma+1}\ordencea N$ a submodel of size
$<\zeta^*$ which contains $|N_{\gamma,1}|\cup |N_{\gamma,2}|$. At limit
stages $\alpha<\zeta$, take
$M_\alpha:=\overline{\bigcup_{\gamma<\alpha}M_{\gamma}}$.

\begin{remark}
Notice that $\langle M_\gamma : \gamma<\zeta \rangle$ is a $\ordencea$-increasing and continuous sequence such that $a\not\hspace{-2.0mm}\indep^{T,\E}_{M_\gamma} M_{\gamma+1}$ for every $\gamma<\zeta$ (because $M_{\gamma+1}$ contains the models that witness the $\E$-tame splitting).
\end{remark}

\indent Let us construct a sequence $\langle M_{\alpha}^*
:\alpha\le \zeta \rangle$ of models and a tree $\langle
h_{\eta}:\eta\in\;^\alpha
2 \rangle$ ($\alpha\le \zeta$) of
$\C{K}$-embeddings such that:

\begin{enumerate}
\item $\gamma<\alpha$ implies $M_\gamma^*\ordencea M_\alpha^*$.
\item $M_{\alpha}^*:=\overline{\bigcup_{\gamma<\alpha}M_{\gamma}^*}$ if
$\alpha$ is limit.
\item $\gamma<\alpha$ and $\eta\in\;^{\alpha}2$ imply that $h_{\eta\rest \gamma}\subset
h_{\eta}$.
\item $h_\eta:M_{\alpha}\to M_{\alpha}^*$ for every $\eta\in \;^\alpha
2$.
\item If $\eta\in\; ^\gamma 2$ then $h_{\eta^\frown 0}(N_{\gamma,1})=h_{\gamma^\frown 1}(N_{\gamma,2})$
\end{enumerate}

Take $M_0^*:=M_0$ and $h_{\langle \rangle}:=id_{M_0}$.
\\ \\
\indent If $\alpha$ is limit, take
$M_{\alpha}^*:=\overline{\bigcup_{\gamma<\alpha}M_{\gamma}^*}$ and
if $\eta\in\hspace{.1mm}^\alpha 2$ define
$h_{\eta}:=\overline{\bigcup_{\gamma<\alpha}h_{\eta\rest \gamma}
}$.
\\ \\
\indent If $\alpha:=\gamma+1$, let
$\eta\in\hspace{.1mm}^{\gamma}2$. Take $\overline{h_{\eta}}\supset
h_\eta$ any automorphism of the monster model $\monster$ (this is
possible because $\monster$ is homogeneous).
\\ \\
\indent Notice that $\overline{h_{\eta}}\circ
F_{\gamma}(N_{\gamma,1})=\overline{h_{\eta}}(N_{\gamma,2})$.
Define $h_{\eta^\frown 0}$ as any extension of
$\overline{h_{\eta}}\circ F_\gamma$ to $M_{\gamma+1}$ and
$h_{\eta^\frown 1}$ as $\overline{h_{\eta}}\rest M_{\gamma+1}$.
Take $M_{\gamma+1}^*\ordencea N$ as any model with density
character $<\zeta^*$ which contains $h_{\eta^\frown l}(M_{\gamma+1})$
for any $\eta\in \hspace{0.1mm}^\gamma 2$ and $l=0,1$.
\\ \\
\indent Take $H_{\eta}$ an automorphism of $\monster$ which
extends $h_{\eta}$, for every $\eta\le\hspace{0.1mm}^{\zeta} 2$.

\begin{claim}\label{LessE2}
If $\eta\neq \nu \in\hspace{.1mm}^{\zeta} 2$ then
$d(\gatp(H_\eta(\overline{a})/M_\zeta^*),
\gatp(H_\nu(\overline{a})/M_\zeta^*))\ge \varepsilon$.
\end{claim}
\bdem Suppose not, then $d(\gatp(H_\eta(\overline{a})/M_\zeta^*),
\gatp(H_\nu(\overline{a})/M_\zeta^*))< \varepsilon$. Let
$\rho:=\eta\land \nu$. Without loss of generality, suppose that $\rho^\frown 0\le \eta$
and $\rho^\frown 1\le \nu$. Let $\gamma:=length(\rho)$. Since
$h_{\rho^\frown 0}(N_{\gamma,1})=h_{\rho^{\frown
}1}(N_{\gamma,2})\ordencea M_{\zeta}^*$, therefore\linebreak
$d(\gatp(H_\eta(\overline{a})/h_{\rho^\frown 0}(N_{\gamma,1})),
\gatp(H_\nu(\overline{a})/h_{\rho^{\frown}1}(N_{\gamma,2}))<
\varepsilon$. Also

\begin{eqnarray*}
d(\gatp(H_{\nu}^{-1}\circ
H_\eta(\overline{a})/F_{\gamma}(N_{\gamma,1})),
\gatp(\overline{a}/N_{\gamma,2})) &=&\\
d(\gatp(H_\eta(\overline{a})/h_{\rho^\frown 0}(N_{\gamma,1})),
\gatp(H_\nu(\overline{a})/h_{\rho^{\frown}1}(N_{\gamma,2})) &<&
\varepsilon
\end{eqnarray*}

(since $H_{\nu}$ is an isometry, $h_{\rho^{\frown}0}=h_{\rho}\circ
F_{\gamma}$, $\rho< \nu$, $\rho^{\frown}0\le \eta$ and
$\rho^{\frown}1\le \nu$). Since $H_{\nu}^{-1}\circ
H_\eta(\overline{a}) \supset F_{\gamma}$, then
$d(F_{\gamma}(p\rest N_{\gamma,1}),p\rest N_{\gamma,2})<\E$, which
contradicts the choice of $N_{\gamma,1}$, $N_{\gamma,2}$ and
$F_{\gamma}$. 
\edem[Claim \ref{LessE2}]

We have that $dc(M_{\zeta}^*)\le \lambda$ (because $dc(M_{\zeta}^*)\le \zeta^*\cdot \zeta =\max\{\mu^+,\zeta\}\cdot \zeta \le \lambda $). Take $M^*\ordenceac M_{\zeta}^*$ of density character $\lambda$; so by claim \ref{LessE2} we have that $dc(\gaS(M^*))\ge 2^{\zeta}>\lambda$, which contradicts
$\lambda$-$d$-stability.
\edem[Prop. \ref{Locality2}]

\begin{proposition}[Weak stationarity over universal models]\label{stationatity2}
For every $\E>0$ there exists $\delta$ such that for every $N_0\ordencea N_1\ordencea N_2$ and every $a,b$, if $N_1$ is universal over $N_0$, $a,b\indep^{T,\delta}_{N_0} N_2$ and
$$\dist(\gatp(a/N_1),\gatp(b/N_1))<\delta,$$ therefore
$$\dist(\gatp(a/N_2),\gatp(b/N_2))<\E.$$
\end{proposition}
\bdem
Take $\delta:=\delta_\E/3$ (see definition of tameness, \ref{tameness}). Let $N^*\ordencea N_0$ be a model of size $<\zeta^*$ which witnesses $a,b\indep^{T,\delta}_{N_0} N_2$. Let $M^{\circ} \ordencea N_2$ be a model of density character $\mu$. Let $M^*\ordencea N_2$ be a model of density character $<\zeta^*$ which contains $|N^*|\cup|M^\circ|$. Since $N_1$ is universal over $N_0$, so it is $<\zeta^*$-universal over $N^*$. Therefore, there exist a model $M'$ such that $N^*\ordencea M'\ordencea N_1$ and an isomorphism $f:M'\stackrel{f}{\cong}_{N^*} M^*$. Since $N^*$ witnesses that $a,b\indep^{T,\delta}_{N_0} N_2$ and $N^*\ordencea M'\stackrel{f}{\cong }_{N^* }M^* \ordencea N_2$, therefore $$\dist(\gatp(a/M^*),\gatp(f(a)/M^*))<\delta$$ and $$\dist(\gatp(b/M^*),\gatp(f(b)/M^*))<\delta.$$ Also, since $f$ is an isometry, by hypothesis we have that

\begin{eqnarray*}
\dist(\gatp(f(a)/M^*),\gatp(f(b)/M^*)) &=&\dist(\gatp(a/M'),\gatp(b/M'))\\
&\le& \dist(\gatp(a/N_1),\gatp(b/N_1))\\
&<&\delta
\end{eqnarray*}

Therefore:

\begin{eqnarray*}
\dist(\gatp(a/M^\circ),\gatp(b/M^\circ)) &\le& \dist(\gatp(a/M^*),\gatp(b/M^*))\\
                                    &\le& \dist(\gatp(a/M^*),\gatp(f(a)/M^*))\\
                                    && + \dist(\gatp(f(a)/M^*),\gatp(f(b)/M^*))\\
                                    && + \dist(\gatp(f(b)/M^*),\gatp(b/M^*))\\
                                    &<& 3\delta=\delta_\E\\
\end{eqnarray*}
By $\mu$-$\mathbf{d}$-tameness, we have that $\dist(\gatp(a/N_2),\gatp(b/N_2))<\E$. \\
\edem[Prop. \ref{stationatity2}] 

\section{A stability transfer theorem}
First, we provide a general stability transfer theorem.
\begin{theorem}\label{stab_transfer2}
Let $\C{K}$ be a 
$\mu$-$\mathbf{d}$-tame (for some $\mu<\kappa$) MAEC. Suppose that 
$\C{K}$ is $[LS(\C{K}),\kappa)$-cofinally-d-stable.
Define $\lambda:=\min\{\theta<\kappa: \mu< \theta \text{ and  $\C{K}$ is $\theta$-$\mathbf{d}$-stable }\}$,  $\zeta:=\min\{\xi: 2^{\xi}>\lambda\}$ and $\zeta^*:=\max\{\mu^+,\zeta\}$. If $cf(\kappa)\ge \zeta^*$ then $\C{K}$ is $\kappa$-$\mathbf{d}$-stable.
\end{theorem}
\bdem
Suppose that this proposition is false. Let $M\in \C{K}$ be a model of density character $\kappa$ such that there are $a_i$ ($i<\kappa^+$) such that\linebreak
$\dist(\gatp(a_i/M),\gatp(a_j)/M)\ge \E$ for every $i<j<\kappa^+$ and for some fixed $\E>0$. {\it Without loss of generality}, we can assume that $M$ is the completion of the union of a $\ordencea$-increasing sequence $(M_i:i<cf(\kappa))$ such that $LS(\C{K})\le dc(M_i)<\kappa$ and $M_{i+1}$ is universal over $M_i$ (this is possible by fact \ref{Existence_Universal} and cofinal-d-stability), for every $i<cf(\kappa)$. By proposition \ref{Locality2}, for every $\E>0$ and every $i<\kappa^+$ there exists $M_{i,\E}\ordencea M$ of density character $<\zeta^*$ such that $a_i\indep^{T,\E}_{M_{i,\E}} M$. Since $dc(M_{i,\E})<\zeta^*\le cf(\kappa)$, there exists $j_i<cf(\kappa)$ such that $M_{i,\E}\ordencea M_{j_i}$. By monotonicity of $\indep^{T,\E}$, we have that $a_i\indep^{T,\E}_{M_{j_i}} M$. By pigeon-hole principle, there exists $i^*<cf(\kappa)$ and $X\subset \kappa^+$ of size $\kappa^+$ such that for every $k\in X$ we have that $a_k\indep^{T,\E}_{M_{j_{i^*}}} M$. By proposition \ref{stationatity2}, there exists $\delta>0$ such that  $\dist(\gatp(a_k/M_{j_{i^*}+1}),\gatp(a_j/M_{j_{i^*}+1}))\ge \delta$. By hypothesis $\C{K}$ is $[LS(\C{K}),\kappa)$-cofinally-d-stable, hence there exists $dc(M_{j_{i^*}+1})\le \theta' <\kappa$ such that $\C{K}$ is $\theta'$-$\mathbf{d}$-stable; we can take $M^*\ordenceac M_{j_{i^*}+1}$ with density character $\theta'$,  so $\dist(\gatp(a_k/M^*),\gatp(a_j/M^*))\ge \delta$ for every $j\neq k\in X$ (this contradicts $\theta'$-$\mathbf{d}$-stability).
\edem[Prop. \ref{stab_transfer2}]

The following corollary lets us go up from d-stability in $\aleph_0$ and $\aleph_1$ to d-stability in $\aleph_n$ for every $n<\omega$.

\begin{corollary}\label{stab_spectrum1}
Let $\C{K}$ be an 
$\aleph_0$-$\mathbf{d}$-tame MAEC. Suppose that $\C{K}$ is $\aleph_0$-$d$-stable and $\aleph_1$-d-stable. 
Then $\C{K}$ is $\aleph_n$-d-stable for all $n<\omega$
\end{corollary}
\bdem
Consider $\mu:=\aleph_0$ and $\kappa:=\aleph_2$. Notice that $\lambda:=\min\{\theta<\kappa: \mu< \theta \text{ and  $\C{K}$ is $\theta$-$\mathbf{d}$-stable }\}=\aleph_1$ and $\zeta:=\min\{\xi: 2^{\xi}>\lambda\}\le \aleph_1$. So, $\zeta^*:=\max\{\mu^+,\zeta\}=\aleph_1$ (independently if CH holds). In this case, $a\indep^T_N M$ (based on $<\zeta^*$-$\E$-non splitting) means that given $\E$ there exists a separable model $N_\E\ordencea N$ such that $a\indep^\E_{N_\E} M$. Notice that $cf(\kappa)=\aleph_2\ge \zeta^*=\aleph_1$, so by theorem \ref{stab_transfer2} we have that $\C{K}$ is $\aleph_2$-d-stable. By an inductive argument, we have that $\C{K}$ is $\aleph_n$-d-stable for all $n<\omega$.
\edem[Cor. \ref{stab_spectrum1}]

The following corollary says that, under the superstability-like assumption below, we can get $\aleph_\omega$-d-stability from d-stability in $\aleph_n$ for every $n<\omega$.

\begin{assumption}[$\E$-locality]\label{superstability_tameness}
For every tuple $\tuple{a}$, every $\E>0$ and every increasing and continuous $\ordencea$-chain of models $\langle M_i : i<\sigma \rangle$, there exists $j<\sigma$ such that $\tuple{a}\indep^{T,\E}_{M_j} \overline{\bigcup_{i<\sigma} M_i}$.
\end{assumption}

\begin{corollary}\label{superstability_spectrum1}
Let $\C{K}$ be a $\aleph_0$-d-tame, $\aleph_0$-d-stable and $\aleph_1$-d-stable MAEC which satisfies assumption \ref{superstability_tameness}. Then $\C{K}$ is $\aleph_\omega$-d-stable.
\end{corollary}
\bdem
By corollary \ref{stab_spectrum1}, $\C{K}$ is $\aleph_n$-d-stable for all $n<\omega$. By reductio ad absurdum, suppose $\C{M}$ is not $\aleph_\omega$-d-stable. So, there exists $M\in \C{K}$ of density character $\aleph_\omega$ such that $dc(\gaS(M))\ge \aleph_{\omega+1}$. {\it Without loss of generality}, we may assume $M$ is the completion of the union of a $\C{K}$-increasing and continuous chain $\{M_n:i<\omega\}$ where $dc(M_n)=\aleph_n$ and $M_{n+1}$ is universal over $M_n$ for all $n<\omega$ (this is possible by fact \ref{Existence_Universal} and $\aleph_n$-d-stability).  So, 
there exist $\E>0$ and $a_i\in \monster$ ($i<\aleph_{\omega+1}$) such that $d(\gatp(a_i/M),\gatp(a_j/M))\ge \E$ for all $i\neq j<\aleph_{\omega+1}$ (we can find them using the same argument when the space is not separable, because $cf(\aleph_{\omega+1})>\omega$, see  \cite{Lima,Wilansky}).
\\ \\
By $\aleph_0$-d-tameness, there exists $\delta_\E>0$ such that for every $p,q\in \gaS(M)$, if $d(p,q)\ge \E$ then there exists $M'\ordencea M$ of density character $\aleph_0$ such that $d(p\rest M',q\rest M')\ge \delta_\E$ (see definition \ref{tameness}). Define $\delta:=\delta_\E/3$.
\\ \\
On the other hand, given $i<\aleph_{\omega+1}$, by the superstability-like assumption \ref{superstability_tameness} there exists $n_i<\omega$ such that $a_i\indep^{T,\delta}_{M_{n_i}} M$. Since $cf(\aleph_{\omega+1})=\aleph_{\omega+1}>\omega$, by pigeon-hole principle there exists a fixed $n<\omega$ and $X\subset\aleph_{\omega+1}$ of size $\aleph_{\omega+1}$ such that $a_i\indep^{T,\delta}_{M_n} M$ for all $i\in X$.
\\ \\
Notice that for every $i\neq j\in X$, $d(\gatp(a_i/M),\gatp(a_j/M))\ge \E$ and\linebreak $a_i,a_j\indep^{T,\delta}_{M_n} M$. We may say that
$$
d(\gatp(a_i/M_{n+1}),\gatp(a_j/M_{n+1}))\ge \delta.
$$
If not, suppose $d(\gatp(a_i/M_{n+1}),\gatp(a_j/M_{n+1}))< \delta$. Let $N^*\ordencea M_{n}$ be a model of size $\aleph_0$ which witnesses $a_i,a_j\indep^{T,\delta}_{M_n} M$. Let $M^{\circ} \ordencea M$ be any model of density character $\aleph_0$. Let $M^*\ordencea M$ be a model of density character $\aleph_0$ which contains $|N^*|\cup|M^\circ|$. Since $M_{n+1}$ is universal over $M_n$, so it is universal over $N^*$. Therefore, there exist a model $M'$ such that $N^*\ordencea M'\ordencea M_{n+1}$ and an isomorphism $f:M'\stackrel{f}{\cong}_{N^*} M^*$. Since $N^*$ witnesses that $a_i,a_j\indep^{T,\delta}_{M_n} M$ and $N^*\ordencea M'\stackrel{f}{\cong }_{N^* }M^* \ordencea M$, therefore $$\dist(\gatp(a_i/M^*),\gatp(f(a_i)/M^*))<\delta$$ and $$\dist(\gatp(a_j/M^*),\gatp(f(a_j)/M^*))<\delta$$
Since $M'\ordencea M_{n+1}$, we have that
\begin{eqnarray*}
\dist(\gatp(a_i/M'),\gatp(a_j/M')) &\le& \dist(\gatp(a_i/M_{n+1}),\gatp(a_j/M_{n+1}))\\
&<&\delta
\end{eqnarray*}
so, 

\begin{eqnarray*}
\dist(\gatp(f(a_i)/M^*),\gatp(f(a_j)/M^*)) &=& \dist(\gatp(a_i/M'),\gatp(a_j/M'))\\
&<& \delta
\end{eqnarray*}

Therefore:
\begin{eqnarray*}
\dist(\gatp(a_i/M^\circ),\gatp(a_j/M^\circ)) &\le& \dist(\gatp(a_i/M^*),\gatp(a_j/M^*))\\
                                    &\le& \dist(\gatp(a_i/M^*),\gatp(f(a_i)/M^*))\\
                                    && + \dist(\gatp(f(a_i)/M^*),\gatp(f(a_j)/M^*))\\
                                    && + \dist(\gatp(f(a_j)/M^*),\gatp(a_j/M^*))\\
                                    &<& 3\delta=\delta_\E\\
\end{eqnarray*}
By $\aleph_0$-$\mathbf{d}$-tameness, we have that $\dist(\gatp(a_i/M),\gatp(a_j/M))<\E$ (contradiction).
\\ \\
Hence $dc(\gaS(M_{n+1}))\ge \aleph_{\omega+1}>\aleph_{n+1}$, contradicting $\aleph_{n+1}$-d-stability. \ \ \ \ \ 
\edem[Cor. \ref{superstability_spectrum1}]

\begin{corollary}[weak superstability]\label{weak_superstability}
Let $\C{K}$ be an $\aleph_0$-d-tame, $\aleph_0$-d-stable and $\aleph_1$-d-stable MAEC, which also satisfies assumption \ref{superstability_tameness} (countable locality of $\E$-splitting). Then $\C{K}$ is $\kappa$-$d$-stable for every cardinality $\kappa$.
\end{corollary}
\bdem
By induction on all cardinalities $\kappa\ge \aleph_0$, we prove that $\C{K}$ is $\kappa$-d-stable.
By hypothesis, we have $\C{K}$ is $\aleph_0$ and $\aleph_1$-d-stable.
\\ \\
Suppose $\C{K}$ is $\lambda$-d-stable for all $\lambda<\kappa$. Notice that $\mu=\aleph_0$,  $\lambda=\min\{\theta>\mu: \C{K} \text{\ is $\theta$-d-stable }\}=\aleph_1$, $\zeta=\min\{\xi: 2^\xi>\lambda\}\le \aleph_1$ and $\zeta^*=\max\{\mu^+,\zeta\}=\aleph_1$. If $cf(\kappa)>\aleph_0$ then $cf(\kappa)\ge \aleph_1=\zeta^*$, then by theorem \ref{stab_transfer2} $\C{K}$ is $\kappa$-d-stable.
\\ \\
If $cf(\kappa)=\omega$, the argument given in corollary \ref{superstability_spectrum1} works for proving that $\C{K}$ is $\kappa$-d-stable. For the sake of completeness, we provide the proof if $cf(\kappa)=\omega$. Let $\Lambda:\aleph_0\to \kappa$ be a cofinal mapping. By hypothesis, $\C{K}$ is $\Lambda(n)$-d-stable. By reductio ad absurdum, suppose $\C{M}$ is not $\kappa$-d-stable. So, there exists $M\in \C{K}$ of density character $\kappa$ such that $dc(\gaS(M))\ge \kappa^+$. Without loss of generality, we may assume $M$ is the completion of the union of a $\ordencea$-increasing and continuous chain $\{M_n:i<\omega\}$ where $dc(M_n)=\Lambda(n)$ and $M_{n+1}$ is universal over $M_n$ for all $n<\omega$ (this is possible by fact~\ref{Existence_Universal} and $\Lambda(n)$-d-stability). Given $\E>0$, let $a_i\in \monster$ ($i<\kappa^+$) be such that $d(\gatp(a_i/M),\gatp(a_j/M))\ge \E$ for all $i\neq j<\kappa^+$. Let $\delta:=\delta_\E/3$ (where $\delta_\E$ is given in definition \ref{tameness} -tameness-). 
On the other hand, given $i<\kappa^+$, by the superstability-like assumption \ref{superstability_tameness} there exists $n_i<\omega$ such that $a_i\indep^{T,\delta}_{M_{n_i}} M$. Since $cf(\kappa^+)=\kappa^+>\omega$, by the pigeon-hole principle there exists a fixed $n<\omega$ and $X\subset\kappa^+$ of size $\kappa^+$ such that $a_i\indep^{T,\delta}_{M_n} M$ for all $i\in X$.
\\ \\
Notice that for every $i\neq j\in X$, $d(\gatp(a_i/M),\gatp(a_j/M))\ge \E$ and\linebreak $a_i,a_j\indep^{T,\delta}_{M_n} M$. So, by the argument given in corollary \ref{superstability_spectrum1} we may say
$$
d(\gatp(a_i/M_{n+1}),\gatp(a_j/M_{n+1}))\ge \delta.
$$
Hence $dc(\gaS(M_{n+1}))\ge \kappa^+>\Lambda(n+1)$, which contradicts $\Lambda(n+1)$-d-stability.
\edem[Cor. \ref{weak_superstability}]

\end{document}